\documentclass{article}
\usepackage{amssymb}
\usepackage{amsmath}
\newtheorem{thm}{Theorem}[section]
\newtheorem{lem}[thm]{Lemma}

\newtheorem{prop}[thm]{Proposition}

\newtheorem{mthm}[thm]{Main Theorem}
\newtheorem{rem}[thm]{Remark}
\newtheorem{defn}[thm]{Definition}

\renewcommand{\arraystretch}{1.2}

\def\B{{\mathbb {B}}}
\def\C{{\mathbb {C}}}

\def\bH{{\mathbb{H}_g}}

\def\Q{{\mathbb Q}}
\def\R{{\mathbb R}}

\def\Z{{\mathbb {Z}}}
\def\Sn{{Spec(\Z[\zeta_n,1/n])}}
\def\Agn{{A_g(n)}}
\def\Agm{{A_g(m)}}
\def\Acgn{{A^*_g(n)}}
\def\Acgm{{A^*_g(m)}}
\def\Aogn{{A^0_g(n)}}

\def\Aocgn{{A^{0*}_g(n)}}
\def\Aocgm{{A^{0*}_g(m)}}

\def\Bgn{{B_g(n)}}
\def\Bcgn{{B^{*}_g(n)}}

\def\Bocgn{{B^{0*}_g(n)}}

\def\abcd{\bigl(\begin{array}{cc} a&b\\c&d \end{array}\bigr)}

\def\AB0D{\bigl(\begin{array}{cc} A&B\\0&D \end{array}\bigr)}
\def\A-BBA{\bigl(\begin{array}{cc} A&-B\\B&A \end{array}\bigr)}

\def\a{{\alpha}}

\def\l{{\lambda}}
\def\m{{\mu}}

\def\o{{\omega}}
\def\O{{\Omega}}

\def\z{{\zeta}}

\title{{On a theorem of Ihara}}
\author{Arash Rastegar}
\begin{document}
\maketitle
\begin{abstract}
Let $p$ be a prime number and let $n$ be a positive integer prime
to $p$. By an Ihara-result, we mean existence of an injection with
torsion-free cokernel from a full lattice in the space of $p$-old
modular forms, into a full lattice in the space of all modular
forms of level $np$. In this paper, we prove Ihara-results for
genus two Siegel modular forms, Siegel-Jacobi forms and for
Hilbert modular forms. Ihara did the genus one case of elliptic
modular forms [Ih]. We propose a geometric formulation for the
notion of $p$-old Siegel modular forms of genus two using
clarifying comments by R. Schmidt [Sch] and then follow
suggestions in an earlier paper [Ra1] on how to prove Ihara
results. We use the main theorem in [Ra1] where we have extended
an argument by G. Pappas to prove torsion-freeness of certain
cokernel using the density of Hecke-orbits in the moduli space of
principally polarized abelian varieties and in the
Hilbert-Blumenthal moduli space which was proved by C. Chai [Ch].
\end{abstract}


\section*{Introduction}

In this paper, we are interested in proving arithmetic results
for modular forms by introducing geometric interpretations of
some analytic aspects of the theory of modular forms. More
precisely, we suggest a geometric definition for the analytic
notion of old Siegel modular forms and then generalize an
arithmetic result of Ihara to the case of genus two Siegel modular
forms, Siegel-Jacobi forms and Hilbert modular forms.

The theorem of Ihara is used by Ribet to get congruences between
elliptic modular forms. The fundamental paper [Ri1] together with
his ideas in [Ri2] which where generalized by himself and many
others lead to an almost complete classification of congruences
between modular forms of different weights and levels in the
genus one case. The result of Ihara is crucial in the procedure
of raising the level.

To generalize this result to the Siegel case, one has to
introduce a geometric characterization for the space of $p$-old
Siegel modular forms. We introduce explicit correspondences, which
we call Atkin-Lehner correspondences, on appropriate Siegel moduli
spaces which generate the whole $p$-old part out of the pull back
copy of modular forms of level $n$ inside those of level $np$ for
square-free integer $n$. In fact, using the Atkin-Lehner
correspondences, we define a map from four copies of the space of
forms of level $n$ to the $p$-old part of forms of level $np$
which turns out to be injective and generate the whole $p$-old
part. By an Ihara-result we mean cokernel torsion-freeness of the
map induced on the specified full lattices in these vector spaces.
This is what Ihara proved in the elliptic modular case [Ih].
Injection of this map is an automorphic fact. But cokernel
torsion-freeness is proved by getting an injection result in
finite characteristic. We generalize a result of G. Pappas to get
this injection using density of Hecke orbits. This density result
is proved by C. Chai [Ch]. The precise statement of our main
result is as follows.

\begin{mthm} Let $p$ be a prime which does not divide the
square-free integer $n$. Atkin-Lehner correspondences induce an
injection from a full lattice in the space of $p$-old Siegel
modular forms, into a full lattice in the space of all modular
forms of level $np$. The cokernel of this map is free of
$l$-torsions for all primes $l$ not dividing
$2np[\Gamma_0(p):\Gamma'(p)]$ with $l-1>k$ where $\Gamma_0(p)$
and $\Gamma'(p)$ are certain congruence subgroups of $Sp(4,\Z)$.
\end{mthm}

This is an implication of theorems 2.2 and 2.4. Using the same
ideas one can also prove an Ihara result for Siegel-Jacobi forms
of genus two. Let ${ B}^{0*}_g(n)$ denotes the compactification
of the universal abelian variety over the Siegel space ${
A}^0_g(n)$ (Find precise definitions in sections 1 and 3).

\begin{thm} Let $p$ be a prime which does not divide $n$. The
Atkin-Lehner correspondences induce a cokernel torsion-free
injection
$$
H^0( { B}^{0*}_2(n)/\Z_l,\o^{\otimes k}\otimes L^{\otimes
m})^{\oplus 4} \to H^0( { B}^{0*}_2(np)/\Z_l,\o^{\otimes
k}\otimes L^{\otimes m})
$$
for all primes $l$ not dividing $2pn[\Gamma_0(p):\Gamma'(p)]$
with $l-1>k$.
\end{thm}

Theorem 3.1 is a more precise statement. To prove an Ihara result
in the case of Hilbert modular forms is much easier, since we only
have two copies of Hilbert modular forms of level $n$ inside
$p$-old Hilbert modular forms of level $np$. Since our ideas work
for the higher dimensional case of Hilbert-Blumenthal moduli
space, we include an Ihara result for Hilbert modular forms as
well. Although, one can easily prove such a result using
appropriate Shimura curves.

\subsection*{Table of contents}
1. The arithmetic and geometry of Siegel spaces
\\ 2. Ihara result for Siegel modular forms
\\ 3. Ihara result for Siegel-Jacobi modular forms
\\ 4. Ihara result for Hilbert modular forms
\\ Appendix
\\ References

\section{The arithmetic and geometry of Siegel spaces}

In this paper, we only work with geometric formulation of Siegel
modular forms. For automorphic motivations on how to define the
space of $p$-old Siegel modular forms we refer to [Ra1].

By the Siegel upper half-space which we denote by $\bH$ we mean
the set of complex symmetric $g$ by $g$ matrices $\O$ with
positive-definite imaginary part. The quotient of Siegel upper
half-space $\bH$ by $Sp(2g,\Z)$ acting on $\bH$ via M\"obius
transformations is a complex analytic stack $A_{g}$ which could
be thought of as the moduli space of principally polarized abelian
varieties. The universal family of abelian varieties over $\bH$
is given by $ A(\O)=\C^g/(\Z^g \oplus \O.\Z^g)$.

A Siegel modular form of weight $k$ is a holomorphic section of
the line bundle ${\mathcal{O}}_{\bH} \otimes
(\wedge^g\C^g)^{\otimes k}$ which we denote of $\o$. When written
formally, this becomes an expression of the form
$f(\O)(dz_1\wedge...\wedge dz_g)^{\otimes k}$ where $f$ is an
$Sp(2g,\Z)$-invariant complex holomorphic function on $\bH$ which
is also holomorphic at $\infty$. For genus $\geq 2$ the latter
condition is automatically satisfied by Koecher principle.

A discrete subgroup $\Gamma$ of $Sp(2g,\Z)$ is called a congruence
subgroup, if it contains $\Gamma(n)$ for some positive integer
$n$ where
$$
\Gamma(n)=\{\gamma \in Sp(2g,\Z)| \gamma \equiv
\bigl(\begin{array}{cc} {I_g}&0\\0&{I_g}
\end{array}\bigr) \hspace{.2in}(mod \hspace{.05in} n) \}.
$$
Here is an example of a congruence subgroup
$$
\Gamma_0(n)=\{\gamma \in Sp(2g,\Z)| \gamma \equiv
\bigl(\begin{array}{cc} *&* \\ 0&* \end{array}\bigr)
\hspace{.2in}(mod\hspace{.05in} n) \}.
$$

The quotients of $\bH$ by congruence subgroups of $Sp(2g,\Z)$ are
called Siegel spaces. They can be thought of as the moduli space
of principally polarized abelian varieties equipped with certain
level structure. Siegel modular forms can be considered as certain
differential forms on Siegel spaces.

The space of Siegel modular forms can also be formulated in the
language of schemes. Let $S$ be a base scheme. A modular form $f$
of weight $k$ is a rule which assigns to each principally
polarized abelian variety $(A/S,\l)$ a section $f(A/S,\l)$ of
$\o_{A/S}^{\otimes k}$ over $S$ depending only on the isomorphism
class of $(A/S,\l)$ commuting with arbitrary base change. Here
$\o_{A/S}$ is the top wedge of tangent bundle at origin of $A$
over $S$.

To define Siegel modular forms of higher level, one should equip
principally polarized abelian varieties with level structures. Let
$\z_n$ denote an $n$-th root of unity where $n\geq 3$. On a
principally polarized abelian scheme $(A,\l)$ over $\Sn$ of
relative dimension $g$ we define a symplectic principal level-$n$
structure to be a symplectic isomorphism
$\a:A[n]\to(\Z/n\Z)^{2g}$ where $(\Z/n\Z)^{2g}$ is equipped with
the standard non-degenerate skew-symmetric pairing
$$
<,>:(\Z/n\Z)^{2g} \times (\Z/n\Z)^{2g} \to \Z/n\Z
$$
$$
<(u,v) , (z,w)> \mapsto u.w^t-v.z^t
$$

Let $S$ be a scheme over $\Sn$. The moduli scheme classifying the
principally polarized abelian schemes over $S$ together with a
symplectic principal level-$n$ structure is a scheme over $S$
and  will be denoted by $\Agn$. The moduli scheme $\Agn$ over $S$
can be constructed from $\Agn$ over $\Sn$ by base change.

$Sp(2g,\Z/n\Z)$ acts as a group of symmetries on $\Agn$ by acting
on level structures. We will recognize these moduli spaces and
their equivariant quotients under the action of subgroups of
$Sp(2g,\Z/n\Z)$ as Siegel spaces. We restrict our attention to
Siegel spaces over $\Sn$. $\Agn$ is connected and smooth over
$\Sn$. The condition $n \geq 3$ is to guarantee that we get a
moduli scheme, instead of getting only a moduli stack. The
natural morphism $\Agn \to \Agm$ where $m$,$n$ are positive
integers $\geq 3$ and $m|n$ is a finite and etale morphism over
$\Sn$.

Let $\Bgn$ denote the universal abelian variety over $\Agn$. The
Hodge bundle $\o$ is defined to be the pull back via the zero
section $i_0:\Agn \to \Bgn$ of the line bundle
$\wedge^{top}\Omega_{\Bgn/\Agn}$. The Hodge bundle is an ample
invertible sheaf on $\Agn$ and can be naturally extended to a
bundle $\o$ on $\Acgn$. We could define the minimal
compactification $\Acgn$ by the formula
$$
\Acgn =proj(\oplus_{k \geq 0}H^0(\Acgn ,\o^{\otimes k})).
$$
The graded ring above is regarded as a $\Z[\z_n,1/n]$-algebra. The
scheme $\Acgn$ is equipped with a stratification by locally
closed subschemes which are geometrically normal and flat over
$\Sn$. Each of these strata is canonically isomorphic to a moduli
space $A_{i}(n)$ for some $i$ between $0$ and $g$. The map $\Agn
\to \Agm$ can be extended uniquely to $\Acgn \to \Acgm$ for
$m|n$. These maps when restricted to strata, induce the
corresponding natural maps between lower genera Siegel spaces
$A_i(n) \to A_i(m)$. The action $Sp(2g,\Z/n\Z)$ on $\Agn$
naturally extends to an action on the compactified Siegel space
$\Acgn$. This action is compatible with the maps $\Acgn \to
\Acgm$ for $m|n$.

Let $K_0(n)$ denote the subgroup of $Sp(2g,\Z/n\Z)$ fixing the $g$
first $(\Z/n\Z)$-basis elements of $(\Z/n\Z)^{\oplus 2g}$ on
which $Sp(2g,\Z/n\Z)$ acts. Since $\Acgn$ is a projective scheme,
we can define the quotient projective schemes $\Aocgn$ to be the
geometric quotient of $\Acgn$ by $K_0(n)$. This quotient provides
us with a compactification of $\Aogn$ which is the moduli scheme
of principally polarized abelian schemes $(A,\l)$ over $\Sn$,
together with $g$ elements in $A[n]$ generating a symplectic
subgroup isomorphic to $(\Z/n\Z)^g$. Again we have natural maps
$\Aocgn \to \Aocgm$ for $m|n$.

We define the Hodge bundle $\o$ on $\Aocgn$ to be the quotient of
the Hodge bundle $\o$ on $\Acgn$ under the action of the
corresponding subgroup $K_0(n)$ of $Sp(2g,\Z/n\Z)$. This is
possible because the line bundle $\o$ on the space $\Acgn$ is
$Sp(2g,\Z/n\Z)$-linearizable. A Siegel modular form of weight $k$
and full level $n$ is a global section of $\o^k$ on $\Acgn$. Over
the complex numbers, this corresponds to a Siegel modular form of
weight $k$ with respect to $\Gamma(n)$. In this paper, by a
Siegel modular form of weight $k$ and level $n$ we mean a global
section of $\o^k$ on $\Aocgn$. This corresponds to the congruence
subgroup $\Gamma_0(n)$.


\section{Ihara result for Siegel modular forms}

The mod-$p$ Bruhat decomposition implies that, we have the
following decomposition
$$
GSp(4,\Q_p)=\coprod B(\Q_p)w_i\Gamma_0(p),
$$
where $B$ is the Borel subgroup and $w_i$ for $i=1$ to $4$ are
running over the representatives of
$W_{Sp(4,\Z_p)}/W_{\Gamma_0(p)}$. We choose the following
representatives
$$
w_1=id,w_2=\left(\begin{array}{cccc} 1&0&0&0\\0&0&0&1\\0&0&1&0\\
0&-1&0&0 \end{array}\right),w_3=\left(\begin{array}{cccc} 0&0&1&0\\
0&1&0&0\\-1&0&0&0\\0&0&0&1 \end{array}\right),
w_4=\left(\begin{array}{cccc} 0&0&1&0\\0&0&0&1\\-1&0&0&0\\
0&-1&0&0 \end{array}\right).
$$
It is suggested in [Ra1] that because of this picture, there
should be $4$ copies of modular forms of level $n$ inside the
$p$-old part which is analytically characterized by the
conjugates $w_i\Gamma_0(p)w_i^{-1}$ for $i=1$ to $4$. All of these
conjugates lie in $Sp(2g,\Z_p)$. We get the copy associated to
$w_1=id$ simply by pulling back level-$n$ forms to level $np$ via
the projection map between moduli spaces. Following the genus one
case, we can get a second copy using the action of universal
$p$-isogeny $w_p$ which is realized on $\Gamma_0(p)$. Now we have
candidates for two of the copies of level-$n$ forms inside the
$p$-old part.

we are interested in finding $4$ algebraic correspondences on
$A_g^0(np)$ such that the image of the pull back copy under the
action of these correspondences gives us all of the $p$-old
copies. The most natural way to look for these correspondences is
to pull back all the $4$ copies of modular forms of level $n$ to
a congruence subgroup which is of a richer geometric structure.
For example, we can pull back to $\Gamma'(n,p)=\Gamma_0(n)\cap
\Gamma'(p)$ where
$$
\Gamma'(p)=\{\gamma\in Sp(4,\Z)|\gamma\equiv
diag(*,*,*,*)\hspace{.05in} (mod \hspace{.05in} p)\}.
$$
The $4$ specified elements of the Weyl group induce $4$
involutions on the moduli space $A'_g(n,p)$ corresponding to
$\Gamma'(n,p)$. Indeed, the conjugations $w_i \Gamma'(p)w_i^{-1}$
stabilize the congruence group. These involutions have a nice
simple interpretation in terms of the moduli property. We could
also work with $\Gamma'(np)$ and the associated Siegel space
$A'_g(np)$.

The Siegel space associated to $\Gamma'(p)$ is the quotient of
the level-$p$ Siegel space by the subgroup of $Sp(2g,\Z/p\Z)$
which stabilizes all of the $2g$ copies of $\Z/p\Z$ in
$(\Z/p\Z)^{2g}$. The conjugations correspond to symplectic
automorphisms which are well defined on the kind of level
structure we are considering here. It is essential to note that
we can not obtain the $4$ copies of $p$-old modular forms by
applying the above $4$ involutions on the direct pull back of
modular forms via the natural map $\pi:A'_g(n,p)\to A_g(n)$.
Because forms in the pull-back are already invariant under all
$w_i$'s. Instead, we shall apply $w_p$ after pulling these forms
back as far as $A^0_g(np)$ and then pull them back to $A'_g(n,p)$.
The second copy we get in this manner, generates a new copy of
$p$-old forms on $A'_g(n,p)$ by applying involutions by Weyl
elements, to the pull back of the second copy generated by $w_p$
which we can push forward down to $A^0_g(np)$. The forth copy can
be constructed by applying $w_p$ again to the final copy. We
propose the following definition for the space of $p$-old Siegel
forms

\begin{defn}
The space of $p$-old Siegel modular forms of level $np$ is
generated by the images of correspondences $id, w_p, \pi'_*\circ
w_2\circ\pi'^*\circ w_p,w_p\circ\pi'_*\circ w_2\circ\pi'^*\circ
w_p$ acting on the pull-back copy of Siegel modular forms of
level $n$ inside those of level $np$. These are called the
Atkin-Lehner correspondences. The space of new forms is defined
to be the orthogonal complement of the space of $p$-old forms
with respect to the Petterson inner product.
\end{defn}

The careful considerations of R. Schmidt shows that this is a
well-defined notion of old-form for square-free $n$ (Look at
table one in [Sch]). Moreover, one can prove that these
correspondences produce all of the $4$ copies we are expecting
inside the space of Siegel modular forms of level $np$. Here is
the precise statement

\begin{thm} Let $p$ be a prime which does not divide $n$. The
correspondences $id, w_p, \pi'_*\circ w_2\circ\pi'^*\circ w_p$
and $ w_p\circ\pi'_*\circ w_2\circ\pi'^*\circ w_p$ induce a
cokernel torsion-free injection
$$
H^0({ A}^{0*}_2(n)/\Z_l,\o^{\otimes k})^{\oplus 4}\to H^0({
A}^{0*}_2(np)/\Z_l,\o^{\otimes k})
$$
for all primes $l$ not dividing $2np[\Gamma_0(p):\Gamma'(p)]$ with
$l-1>k$.
\end{thm}

This theorem is the main theorem in [Ra1] which is obtained by
generalizing the following result of G. Pappas [Pa]. One shall
note that arguments in [Ra1] work only for genus $2$. The general
case is treated in the new version [Ra2].

\begin{thm}[Pappas] Let $F$ be a field of characteristic zero or
finite characteristic $q$ with $q$ not dividing $pn$. Let $f$ and
$g$ be modular forms on $A_g^0(pn)/F$ which are pulled back from
$A^0_g(n)$. Then $w_p.f+g=0$ implies $f=g=0$ except when $p-1|k$
and the characteristic is finite. In this case, the kernel is
generated by $(H^m,-H^m)$ where $H$ is the Hasse invariant.
\end{thm}

For interested reader, we will reproduce the result of G. Pappas
with a few modifications and simplified proofs in the appendix.

\begin{thm}
Let $n$ be a square-free integer. The space of $p$-old Siegel
modular forms is generated by Siegel modular forms of level $np$
which are eigenforms of Hecke operators and appear as Siegel
eigenforms with the same eigenvalues in the space of Siegel forms
of level $n$.
\end{thm}
\textbf{Proof}. There could be at most four copies of Siegel
modular forms of level $n$ inside $p$-old Siegel modular forms of
level $np$. The four correspondences commute with all Hecke
operators $T_q$ with $q$ prime to $p$. Therefore, by corollary
5.3 of [Sch] they produce old forms. By theorem 2.2 the four
images are linearly independent. Therefore, every $p$-old Siegel
modular form is generated by the images of these correspondences.
$\Box$

Now, it is clear that our main theorem follows from theorems 2.2.
and 2.4.
\section{Ihara result for Siegel-Jacobi modular forms}

The universal abelian variety $\Bgn$ hosts the line bundle
$\wedge^*\Omega_{\Bgn/\Agn}$ which can be naturally extended to a
bundle on the universal semi-abelian variety $\Acgn$. We denote
this line bundle by $\pi^*\o$ where $\pi$ is the projection from
the universal abelian variety to the base. Unlike the Hodge
bundle $\pi^*\o$ is not an ample invertible sheaf on $\Bgn$.

One can construct an ample line bundle on the universal abelian
variety $B_{g,n}$ via the principal polarization of $\Bgn$. Let
$$
\Delta :\Bgn \to \Bgn{\times}_{\Agn}\Bgn
$$
denote the diagonal map. We set $L=\Delta ^*P_{\Bgn}$ where $P$
is chosen to be the invertible sheaf which is twice the principal
polarization of $\Bgn$. The invertible sheaf $\pi^*\o^{\otimes k}
\otimes L^{\otimes m}$ is ample on $\Bgn/\Sn$ for $k,m \gg 0$.

Faltings and Chai construct an arithmetic compactification of the
fiber power of universal abelian variety over the base (theorem
1.1 of chapter VI in [Fa-Ch]). Let $Y=\Bgn^{\times s}_{\Agn}$ and
$\pi :Y \to \Agn$ denote the natural projection to the base and
let $T_G=Lie( G_{g,n})^s$. Then there exists an open embedding $Y
\to \widehat Y$ and a proper extension of the projection over the
base $\hat{\pi}:\widehat Y \to \widehat A_{g}(n)$ such that
\\i)$\widehat Y$ and $\widehat A_{g}(n)$ are smooth over
$\Z[1/n,\z_n]$ and the complement $\widehat Y-Y$ is a relative
divisor with normal crossings.
\\ii)The translation action of $Y$ on itself extends to an action of $
G_{g,n}^s$ on $\widehat Y$.
\\iii)$\widehat{\Omega}^1_{Y/\Agn}= \widehat{\Omega}^1_{Y}/\hat
{\pi}^*(\widehat{\Omega}^1_{\Agn})$ is locally free and
isomorphic to $\hat {\pi}^*(T_{ G_{g,n}^{*s}})$.
\\iv)$R^a \hat {\pi}_*(\wedge^b \widehat{\Omega}^*_{Y\Agn})=
(\wedge^a T_{ G_{g,n}}^s)\otimes(\wedge^b T_{ G_{g,n}}^{*s}).$
\\where all these isomorphisms extend the canonical isomorphisms
over $\Agn$. Here $\widehat{\Omega}^1_{Y/S}$ is an abbreviation
for $\Omega^1_{\widehat{Y}/S}[dlog\infty ]$. In case $s=1$ we
denote this compactification by $\Bcgn$. The line bundles
$\pi^*\o$ and $L$ naturally extend over this compactification.
The compactified universal abelian variety $\Bcgn$ does not
contain the universal semi-abelian variety over the base. This
compactification is $Sp(2g,\Z/n\Z)$-linearizable and one gets
compactifications $\Bocgn$ over $\Aocgn$.

Let $m,n,k \in \mathbb N$ and let $R$ be a $\Z[1/n,\z_n]$-module.
A Siegel-Jacobi form $f$ of genus $g$, weight $k$, index $m$ and
full level $n$ with coefficients in $R$ is an element in
$$
J_{k,m}(n)(R)=H^0( B_g(n)/R,\pi^*\o^{\otimes k}\otimes L^{\otimes
m} ).
$$
The action of double-cosets as geometric correspondences does not
preserve $J_{k,m}(n)$. However the action of the Hecke algebra
$H^{(2)}_p$ associated to the pair
$$
( GSp^{(2)}(2g,\Q_p), GSp^{(2)}(2g,\Z_p))
$$
can be defined [Kr]. One can show that this Hecke algebra is
isomorphic to the algebra generated by connected components of
$p^2$-isogenies. By $p^2$-isogenies, we mean isogenies of over
$\Z[1/p]$ between two principally polarized abelian scheme over
$\Z[1/p]$ which are of degree $P^{2eg}$ where $e$ is a positive
integer. Each connected component of $p^2$-isogenies acts on
$J_{k,m}(n)$ via the following diagram

\begin{equation}
\label{diagram}
\renewcommand{\arraystretch}{1.2}
\arraycolsep=-2pt
\begin{array}{ccccc}
& \pi_1^* B_{g,n} & \buildrel \pi_0 \over \longrightarrow
& \pi_2^* B_{g,n}\\
\pi_1'\swarrow & & \searrow \swarrow & & \searrow\pi_2'\\
 B_{g,n} & & Z & &  B_{g,n} \\
\searrow & & \pi_1\swarrow \searrow\pi_2 & & \swarrow \\
& \Agn & \buildrel id \over \longrightarrow & \Agn \\
\end{array}
\end{equation}
where $\pi_0$ is the universal isogeny over $Z$ and the maps
$\pi_1$ and $\pi_2$ are natural projections from the appropriate
Siegel space $Z$. If $\phi$ is an isogeny of degree $p^{2eg}$,
then we have an isomorphism $ \phi^*\pi_2'^*L\cong
\pi_1'^*L^{\otimes p^{2e}}$. We define the action of $Z$ on
$J_{k,m}(n)$ in the following manner
$$
J_{k,m}(n) \buildrel \pi_2'^* \over \longrightarrow H^0(\pi_2^*
B_{g,n},\pi_2'^*(\o^{\otimes k}\otimes L^{\otimes m})) \buildrel
\Phi^* \over \longrightarrow H^0(\pi_1^*
B_{g,n},\pi_1'^*(\o^{\otimes k}\otimes L^{\otimes p^{2e}.m}))
$$
$$
\buildrel \pi_{1*}' \over \longrightarrow H^0(
B_{g,n},\o^{\otimes k}\otimes L^{\otimes p^{2e}.m}) \cong H^0(
B_{g,n},\o^{\otimes k}\otimes [P^2]^*L^{\otimes m}) \to J_{k,m}(n)
$$
The last mapping is possible to be defined because of finiteness
and flatness of multiplication by $p$.

In order to get appropriate correspondences on the universal
abelian variety $\Bocgn$ we construct compactifications
$B'^*_{g}(n,p)$ of the universal abelian variety $B'_{g}(n,p)$
over the Siegel spaces $A'_{g}(n,p)$. In the particular case of
genus two, the involutions $w_i$ for $i=1$ to $4$ on $A'_{2}(n,p)$
coming from the Weyl elements, induce $4$ involutions on the
compactified universal abelian varieties which will be denoted
again by the same notation $ w_i:B'^*_{2}(n,p) \to B'^*_{2}(n,p)$.
The universal $p^2$-isogeny on $A^{0*}_{2}(n)$ can be extended to
a map $ w_p:B^{0*}_{2}(np) \to B^{0*}_{2}(np)$. We extend the
projection $\pi':A'_{g}(n,p)\to A^0_g(np)$ to $\pi':B^0_{g}(np)\to
B^0_g(np)$.

\begin{thm} Let $p$ be a prime which does not divide $n$.
The geometric correspondences $id, w_p, \pi'_*\circ
w_2\circ\pi'^*\circ w_p$ and $ w_p\circ\pi'_*\circ
w_2\circ\pi'^*\circ w_p$ induce a cokernel torsion-free injection
$$
H^0( { B}^{0*}_2(n)/\Z_l,\o^{\otimes k}\otimes L^{\otimes
m})^{\oplus 4} \to H^0( { B}^{0*}_2(np)/\Z_l,\o^{\otimes
k}\otimes L^{\otimes m})
$$
for all primes $l$ not dividing $2pn[\Gamma_0(p):\Gamma'(p)]$
with $l-1>k$.
\end{thm}
\textbf{Proof}. One can restrict Siegel-Jacobi forms to the base
in order to get Siegel modular forms. This geometric restriction
commutes with the action of correspondences and the map between
cohomology groups. The restricted map with coefficients in a
field $F$
$$
H^0(A^{0*}_2(n)/F,\o^{\otimes k}\otimes L^{\otimes m})^{\oplus 4}
\to H^0(A^{0*}_2(np)/F,\o^{\otimes k}\otimes L^{\otimes m})
$$
is proved to be injection for $F$ of characteristic zero or
characteristic $l$ for $l$ not dividing
$2pn[\Gamma_0(p):\Gamma'(p)]$ with $l-1>k$. This implies
injectivity of the map between Siegel-Jacobi forms with
coefficients in $F$, because by applying multiplication by $l^n$
one can get Siegel-Jacobi forms of higher index restricting to a
nonzero Siegel modular form on the base. This implies cokernel
torsion-freeness of the map between Siegel-Jacobi forms.$\Box$

\section{Ihara result for Hilbert modular forms}

Here, we consider Hilbert modular forms in geometric framework
(look at [Wi] for both geometric and analytic approaches). Fix a
totally real number field $F$ with $g=[F:\mathbb Q ]>1$. Let
$O=O_F$ denote the ring of integers of $F$ and $D_F$ its
different. Also, fix a basis ${a_1,...,a_g}$ for $O_F$ over $\Z$
and thus an isomorphism
$$
O_F/nO_F \to (\Z/n\Z)^g.
$$
We are interested in pairs $(X/s,m)$
consisting of an abelian scheme $X/S$ over $S$ of relative
dimension $g$ and a homomorphism $m:O_F\to End(X/S)$ such that
the relative Lie algebra $Lie(X/S)$ is locally a free
$O_F\otimes_{\Z} O_S$-module of rank one on the Zariski site. Let
$c^{\sharp}=( c, c_+)$ be a projective rank-1 $O_F$-module
equipped with a notion of positivity, i.e. an ordering of the
rank one free $\R$-module $ c\otimes_{\tau ,O_F}\R$ for each
embedding $\tau :F \hookrightarrow \R$. A
$c^{\sharp}$-polarization is an isomorphism of etale sheaves of
$O_F$-modules on $S$ with positivity
$$
\lambda:(P,P_+)\buildrel
\sim \over \longrightarrow ( c, c_+)
$$
where $P=P(X,m)$ is the
sheaf $T/S\mapsto Hom_{T,O}(X_T,X^t_T)_{sym}$ of $O$-linear
quasi-polarizations of $X$ with $P_+$ the subsheaf of $O$-linear
polarizations of $X$.

A principal level-$n$ structure on the pair $(X/s,m)$ is an
$O$-linear symplectic isomorphism
$$
\alpha:X[n]\to (O_F/nO_F)^2=(\Z/n\Z)^{2g},
$$
where $X[n]\subset X$ is the subgroup scheme of the torsion points
of $X$. $\a$ induces an isomorphism
$$
P(X,m)\otimes_{O_F} \wedge^2_{O_F}X[n] \buildrel \sim \over
\longrightarrow D_F^{-1}/nD_F^{-1}(1).
$$
Let $n$ be an integer $\geq 3$. Let $S$ be a
$Spec(\Z[\z_n,1/n])$-scheme. The moduli scheme over
$Spec(\Z[\z_n,1/n])$ which classifies $c^{\sharp}$-polarized pairs
$(X,m)$ over $S$ which are equipped with a principal level-$n$
structure is denoted by $M(c^{\sharp},n)$. This is a smooth and
faithfully flat scheme of relative dimension $g$. Let $m|n$ be an
integer $\geq 3$. The natural morphism
$$
M(c^{\sharp},n)\to
M(c^{\sharp},m)
$$
is finite and flat over $Spec(\Z[\z_n,1/n])$.
The moduli scheme $M(c^{\sharp},n)$ depends on the isomorphism
class of $c^{\sharp}$ in the ideal class group $Cl(F)$ of $O_F$.

The moduli space with polarization varying in the class group is
the full Hilbert-Blumenthal moduli scheme $M(pol,n)$ with
level-$n$ structure. Like the Siegel case, there exists a
universal abelian variety over $M(c^{\sharp},n)$. We define the
Hodge bundle $\o$ as the top wedge of the relative dualizing
sheaf. The Hodge bundle is an an ample line bundle. Koecher
principle will still remain valid. Namely sections of
$\o^{\otimes k}$ over $M(c^{\sharp},n)$ naturally extend to
honest sections on the compactification $M^*(c^{\sharp},n)$.
Indeed, we have a natural compactification
$$
M^*(c^{\sharp},n)=Proj(\sum_{k\geq 0}
H^0(M(c^{\sharp},n),\omega^{\otimes k}))
$$
where the graded ring is considered as an algebra over
$\Z[\z_n,1/n]$. This is a canonical projective scheme
$M^*(c^{\sharp},n)$ of finite type over $Spec(\Z[\z_n,1/n])$,
containing $M(c^{\sharp},n)$ as a dense open subscheme. The
compliment is finite and flat. In fact this complement is
isomorphic to a disjoint union of finitely many copies of
$Spec(\Z[\z_n,1/n])$. A cusp of $M(c^{\sharp},n)$ is
characterized by an extension of projective $O_F$-modules
$$
0\to D_F^{-1} a^{-1}\to K \to  b \to 0
$$
and an $O_F$-linear isomorphism
$$
\gamma:K/nK \buildrel \sim \over
\longrightarrow (O_F/nO_F)^2
$$
where $ a$ and $ b$ are rank-1
projective modules equipped with an isomorphism
$$
\beta:b^{-1}a \buildrel \sim \over \longrightarrow  c.
$$
The number of cusps of $M(c^{\sharp},1)$ is equal to the class
number of $O_F$. From this data, for every cusp of
$M(c^{\sharp},n)$ we get an isomorphism
$$
c\otimes_{\Z} \Z/n\Z
\buildrel \sim \over \longrightarrow D_F^{-1}\otimes_{\Z} \Z/n\Z
$$
induced by $\beta$.

The group $Sp(2,O_F/nO_F)$ acts naturally on $M^*(c^{\sharp},n)$
by acting on the level structure at each point. This action is
compatible with the morphism
$$
M(c^{\sharp},n)\to M(c^{\sharp},m)
$$
for $m|n$. We can consider quotients of the moduli scheme by
subgroups of $Sp(2,O_F/nO_F)$. We can define a compactification
$M^{0*}(c^{\sharp},n)$ for $M^0(c^{\sharp},n)$ by normalization
of
$$
M^0(c^{\sharp},n)\to M^*(c^{\sharp},n).
$$

Let $R$ be a $\Z[\z_n,1/n]$-algebra. We define a Hilbert modular
form of genus $g$ and weight $k$, with coefficients in $R$
regular at cusps, to be an element of
$H^0(M(pol,n),\omega^{\otimes k}\otimes R)$. Koecher principle is
also valid for Hilbert modular forms. In other words, any Hilbert
modular form with coefficients in $R$ extends to a unique section
of $\o^{\otimes k}$ on $M^*(pol,n)$. From the moduli property,
one can define an embedding of $M^*(pol,n)$, the compactified
moduli stack of polarized pairs $(X,m)$ into $A^*_g$ taking cusps
to cusps, such that $\omega$ on $A^*_g$ restricts to $\omega$ on
$M^*(pol,1)$. This way, one can get Hilbert modular forms by
restriction of Siegel modular forms to the Hilbert-Blumenthal
moduli space [Wi].

Dividing by subgroups of $Sp(2,O_F/nO_F)$ acting on $M^*_g(n)$
one gets etale quotients of Hilbert-Blumenthal moduli space of
level $n$. Consider the congruence subgroup of $SL(2,O_F)$ defined
by
$$
\Gamma^0( n)=\{\a =\abcd\in GL_2^+(F)|a,d\in O_F,b\in
{\delta}^{-1}, c\in  {n\delta}, ad-bc\in O_F^{\times}\}.
$$
Note that
$$
\bigl(\begin{array}{cc} 0&-1\\1&0
\end{array}\bigr)\Gamma^0( p)\bigl(\begin{array}{cc} 0&-1\\1&0
\end{array}\bigr)^{-1} = \bigl(\begin{array}{cc} 1&0\\0& {p}^{-1}
{\delta}^{-2}\end{array}\bigr)\Gamma^0( p)
\bigl(\begin{array}{cc} 1&0\\0& {p}^{-1} {\delta}^{-2}
\end{array}\bigr)^{-1}.
$$
Thus, we have an involution $w_{ p}$ on the moduli space
$M^0(pol,{pn})$ associated to $\Gamma^0({pn})$ whose level
structure is defined by fixing a symplectic subgroup of abelian
scheme of $np$-torsion points $(O_F/ {np}O_F)^2$ isomorphic to
the group scheme $(O_F/{p}O_F)$. Having $w_{ p}$ we get a map
$$
\pi^*\oplus w_p\circ\pi^*:H^0(M^0(pol, n),\o^{\otimes
k}/O_F)^{\oplus 2} \to H^0(M^0(pol, {pn}),\o^{\otimes k}/O_F)
$$
\begin{thm} Let $p$ be a prime ideal which does not divide $n$. There
exists a cokernel torsion-free injection
$$
H^0(M^0(pol, n)/O_{F,\ell},\o^{\otimes k})^{\oplus 2} \to
H^0(M^0(pol, {pn})/O_{F,\ell},\o^{\otimes k})
$$
for all primes $l$ of the ring $O_{F}$ not dividing $2pn$ with
$l-1>k$. Here $O_{F,\ell}$ denotes the $\ell$-adic localization of
$O_F$ and $l$ is the characteristic of its residue field.
\end{thm}

The proof is in the same lines as the appendix. But, one can
easily prove such a result using appropriate Shimura curves. This
is why we leave the proof to interested reader.

\section*{Appendix}

We shall thank G. Pappas, who communicated his results to us
before publication. His results works both in finite
characteristics and characteristic zero. The idea is essentially
due to F. Diamond and R. Taylor [Di-Ta]. Pappas uses density of
Hecke orbits to apply the ideas of [Di-Ta] in a much more general
setting. Density of Hecke orbits was proved by C. Chai [Ch].

Many results in this section do not hold, when $g=1$.
We assume $g \geq 2$ through the whole section. Let $F$ be a field of
characteristic 0 or of finite characteristic $q$, not dividing the level $np$.

\begin{prop}
The line bundle $\o^{\otimes k}$ on $A^0_g(n)$ is nontrivial for
nonzero $k$.
\end{prop}
\textbf{Proof}. Moret-Bailly and Oort introduce a principally
polarized abelian scheme $A'\to \mathbb P^1(F)$ of relative
dimension two with restriction of $\o$ on $\mathbb P^1(F)$
isomorphic to $O(q-1)$. The local system of $n$-torsion points is
trivial over $\mathbb P^1(F)$ thus by adding a constant family of
right dimension to $A'$ we can induce an embedding $\mathbb
P^1(F)\to A^0_g(n)$. This implies the claim.$\Box$

\begin{rem}. On $X^0(3)$ we have $\o =O_{\mathbb P^1(F)}(1)$ and on
$X^0(4)$ we have $\o =O_{\mathbb P^1(F)}(3)$. In both cases we can
have a divisor for $\o$ supported in the cusps.$\Box$
\end{rem}

Let $S$ be a scheme of characteristic $q$, and $A\to S$ be a
semi-abelian scheme over $S$. Frobenius morphism induces a
homomorphism
$$
Fr^*:e^*\Omega^1_{A/S}\to (e^*\Omega^1_{A/S})^{(q)}
$$
where $(e^*\Omega^1_{A/S})^{(q)}$ is the same as
$e^*\Omega^1_{A/S}$ except that, the underlying $ O_S$-module
structure is combined with Frobenius. Therefore, we get a
homomorphism $det(Fr^*):\o_{A/S}\to \o^{\otimes q}_{A/S}$. Since
$det(e^*\Omega_{A/S}^{1(q)})=det(e^*\Omega^1_{A/S})^{\otimes q}$,
over $\Z/q\Z$ this homomorphism corresponds to a section $H$ of
$\o^{\otimes (q-1)}$.

\begin{thm} Let $\pi:A^0_g(np) \to A^0_g(n)$ be the projection. Then
the morphism
$$
\pi^*\oplus w_p\circ\pi^*:H^0(A^0_g(n)/F,\o^{\otimes k})\oplus
H^0(A^0_g(n)/F,\o^{\otimes k}) \to H^0(A^0_g(np)/F,\o^{\otimes k})
$$
is injective if $F$ is of characteristic zero, or if $k$ is not
divisible by $q-1$. In case $k\equiv 0\hspace{.05in}(mod
\hspace{.05in} q-1)$ the kernel is generated by $(H^m,H^{-m})$
where $H$ is the Hasse invariant and m=k/(q-1).
\end{thm}
\textbf{Proof}. Injection in characteristic zero can be reduced
to $\C$, and in characteristic $q$ can be reduced to $\overline
{\mathbb F}_q$ by flat base change. Suppose $(f_1,f_2)$ be an
element in the kernel of $\pi^*\oplus w^0_p\circ\pi^*$. Consider
an ordinary point on $A^0_g(np)/F$. It corresponds to an isogeny
$\phi:A\to B$ between principally polarized abelian varieties.
Let $\o_A$ and $\o_B$, be a local bases for $\o_{A/F}$ and
$\o_{B/F}$, respectively. We shall regard $f_1$ and $f_2$ as
forms on $A^0_g(np)$ which are pulled back from $A^0_g(n)$. But
when we address the values of them, we shall not specify the
$p$-isogeny associated to the abelian variety on that point,
because these forms being pull backs are independent of these
isogenies. Then by definition of $w_p$, we have
$f_1(A,\phi^*\o_B)=-f_2(B,\o_B)$, and also
$f_1(B,\phi^{t*}\o_A)=-f_2(A,\o_A)$ for the dual isogeny
$\phi^t:B\to A$. If $\psi: B\to C$ be another isogeny of the same
type, we have $f_1(B,\psi^*\o_C)=-f_2(C,\o_C)$. Therefore, $\l
f_2(A,\o_A)=\m f_2(C,\o_C)$ for non-zero constants constants $\l$
and $\m$. We conclude that for a chain of isogenies $A\to B \to
C$ of degree $P^g$ corresponding to points in $A^0_g(np)$, if
$f_2$ vanishes at a point in $A^0_g(np)$ over $A\in A^0_g(n)$,
then it also vanishes at every point in $A^0_g(np)$ which is over
$C\in A^0_g(n)$. Let $Z$ denote the subscheme of $A^0_g(n)$
defined by $f_2=0$. We shall show in the following lemmas that
$Z=A^0_g(n)$. This will prove that $f_1=f_2=0$ and we get
injection.

In case $f_2$ only vanishes at supersingular points, there will
be a positive integer $m$ such that, $f_2/H^m$ is a global
non-zero section of $\o^{m-k(q-1)}$ on $A^0_g(n)$. This implies
triviality of the bundle $\o^{m-k(q-1)}$ and thus $k=m(q-1)$. Let
$\overline {S}_g(n)$ denote the divisor of $H$ in $\overline
{A}^0_g(n)$ where $\overline {A}^0_g(n)$ denotes a smooth
toroidal compactification of $ {A}^0_g(n)$ (Chapter IV in
[Fa-Ch]) and let $S_g(n)$ denote the divisor of $H$ in $A^0_g(n)$.
$\overline {S}_g(n)$ is the reduced Zariski closure of $S_g(n)$.
which is of pure codimension one in $\overline {A}^0_g(n)$, by
construction of toroidal compactifications. The complement
$A^0_g(n)-S_g(n)$ is the ordinary locus. So in this situation, we
get $f_2=\l .H^m$ for some $\l \in F$.$\Box$

\begin{lem}
The subscheme $Z$ of $A^0_g(n)$ satisfies the following condition:
\\(*)If $\phi :A\to B$ is a principally polarized $p$-isogeny
over $F$ of degree $p^c$ with $2g|c$ and one of the points
corresponding to $A$ or $B$ belong to $Z(F)$, then so does the
other point.
\end{lem}
\textbf{Proof}. [Sh] proves that for a pricipally polarized
isogeny $\phi: A \to B$ of degree $p^{cg}$ over $F$ there exists
a chain of principally polarized isogenies of degree $p^g$ whose
composition equals to $\phi$. The details are left to the
reader.$\Box$

\begin{lem}
There exists a closed subscheme $Z'$ of $A^0_g(n)$ which contains
$Z$ and satisfies the following property
\\(*')If $\phi :A\to B$ is a principally polarized $p$-isogeny over
$F$ and one of the points corresponding to $A$ or $B$ belong to
$Z'(F)$ then so does the other point.
\end{lem}

\begin{rem}. By results of [Ch] any closed subscheme of $A^0_g(n)$
satisfying the above property is dense in $A^0_g(n)$.$\Box$
\end{rem}
\textbf{Proof}. We are considering the case that $Z$ has an
ordinary point. Suppose $Z$ be of codimension one. Let $Z''$
denote the Zariski closure of the set of all ordinary points
which are connected to a point on $Z$ by an isogeny of degree
dividing $p^{2g}$. For every point in $Z$ there are at most
finitely many points connected to it by isogenies. Thus $Z''$ has
also codimension one in $A^0_g(n)$. Now define $Z'=Z''\cup
S_g(n)$. $Z'$ satisfies the property mentioned above.$\Box$

\begin{lem} Let $f$ and $g$ be forms as above of the same weight.
To show that $w_p.f+g=0$ implies $f=g=0$ it is enough to prove
$f=w_p.f$ implies $f=0$ for arbitrary weight.
\end{lem}
\textbf{Proof}. This is because $w_p$ is an involution.
$w_p.f+g=0$ implies that $(f+g)^{\otimes 2}$ and $(f-g)^{\otimes
2}$ are both $w_p$-invariant. If being $w_p$-invariant implies
vanishing, we get $(f+g)^{\otimes 2}=(f-g)^{\otimes 2}=0$ and
thus $f+g=f-g=0$ which implies vanishing of $f$ and $g$.$\Box$
\\ \textbf{Proof(2.3)}. Let $h$ be a form on $A^0_g(n)$ whose pull
back to $A^0_g(pn)$ is $w_p$-invariant. Note that if $A\to B$ is
an isotropic isogeny of degree $p^g$ between two ordinary
principally polarized abelian varieties then if $h$ vanishes at
the point corresponding to $A$ it also vanishes at the point
corresponding to $B$. So by the argument of Pappas, the subscheme
$Z$ defined by $h=0$ contains the full Hecke orbit of
 each of its ordinary points. By a result of [Ch] we know that each Hecke
orbit is dense. Hence $Z=A_g$ and $h=0$. Except if $Z$ does not
contain any ordinary point. In this case we repeat the argument
of Pappas and deduce that $h=\l .H^m$ and $p-1|k$.$\Box$

\subsection*{Acknowledgements}

We have benefited from conversations with G. Pappas, A. Rajaei,
C. Skinner, R. Takloo-bighash and A. Wiles. We shall thank them
heartily. We also thank Princeton University, Institute for
Studies in Theoretical Physics and Mathematics (IPM) and Sharif
University of Technology for their partial support. Particular
thanks go to a referee who introduced [Sch] to the author.


\thispagestyle{plain}        
\addcontentsline{toc}{chapter}{Bibliography}
\bibliographystyle{amsplain}
\bibliography{Bibliography}

Sharif University of Technology, e-mail: rastegar@sharif.edu
\\Institut des Hautes Etudes Scientifiques, e-mail:
rastegar@ihes.fr

\end{document}